\documentclass[10pt]{article}

      \usepackage{latexsym}
         \usepackage[reqno, namelimits, sumlimits]{amsmath}
         \usepackage{amssymb, amsfonts}

         \usepackage{amsthm}

\newtheorem{theorem}{Theorem}[section]
\newtheorem{lemma}[theorem]{Lemma}
\newtheorem{remark}[theorem]{Remark}
\newtheorem{proposition}[theorem]{Proposition}

\DeclareMathOperator{\ep}{\varepsilon}
\DeclareMathOperator{\OOa}{\mathcal O_\theta}
\DeclareMathOperator{\OOd}{\mathcal O_{\delta}}

\newcommand{\be}{\begin{equation}}
\newcommand{\ee}{\end{equation}}

\newcommand{\rf}[1]{(\ref{#1})}

\newcommand{\Om}{\Omega}

\newcommand{\vf}{\ensuremath{\varphi}}

\renewcommand{\Re}{\operatorname{Re}}

\newcommand{\R}{{\mathbb R}}

\newcommand{\OO}{\mathcal O}

\newcommand{\dd}{\Delta}

\newcommand{\nn}{\nabla}

\newcommand{\homog}[1]{$m-$homogeneous}

\newcommand{\oth}{\OO_\theta}

\title{
Backward uniqueness for the heat equation in cones}

\author{Lu Li\footnote{University of Minnesota; supported in part by Louise T. Dosdall fellowship and Doctoral Dissertation Fellowship}
\and Vladim\'{\i}r \v Sver\'ak\footnote{University of Minnesota; supported in part by NSF grant DMS 0800908} }
\date{}
\begin{document}

\maketitle

\begin{abstract} It was shown in \cite{ESS, SS1} that a bounded solution of the heat equation in a half-space which becomes zero at some time must be identically zero, even though no assumptions are
made on the boundary values of the solutions. In a recent example, Luis Escauriaza showed that
this statement fails if the half-space is replaced by cones with opening angle smaller than $ 90^\circ$.
Here we show that the result remains true for cones with opening angle larger than $110^\circ$.
\end{abstract}

\section{Introduction}

Consider an open set $\Om\subset \R^n$.
Let $u$ be a bounded solution of the equation
\be\label{heatIntro}
u_t-\dd u + b(x,t)\nn u + c(x,t) u = 0 \quad\mbox{ in $\Om\times(0,T)$},
\ee
 where the coefficients $b=(b_1,\dots, b_n),\,c$ are measurable and bounded. 
  We say that $\Om$ has the {\it backward uniqueness property}
if the following statement holds:

\smallskip

\noindent
{(BU)}$\,\,\,$ {\sl If a bounded $u\colon\Om\times(0,T)\to\R$ satisfies~\rf{heatIntro} and $u(\,\cdot\,, T)=0$, then $u\equiv 0$ in $\Om\times(0,T)$.}

\smallskip
\noindent

 It is important to emphasize that
 no assumptions are made about $u$ at the parabolic  boundary
$\partial \Om\times(0,T)\cup (\Om\times\{0\})$. In fact, we can think about the problem
in terms of control theory: we are given some initial data $u_0\colon\Om\to \R$, and we wish
to find a suitable boundary condition $g$ on the lateral boundary $\partial \Om\times(0,T)$
so that when we solve equation~\rf{heatIntro} with $u_0$ as the initial condition and $g$ as the boundary
 condition, the solution will become exactly zero at time $t=T$. In this interpretation condition
 (BU) means that we can never achieve the exact boundary control of any non-trivial initial condition.
 
While control theory for PDEs seems to be the most natural background for (BU),
the problem also appeared in regularity theory
of parabolic equations,  such as the Navier-Stokes equations, harmonic map
heat flows, or semi-linear heat equations, see~\cite{ESS, Mizoguchi, Wang}. 
The specific unbounded domains which arise in this connection are
complements of closed balls (for interior regularity), half-spaces (for boundary
regularity at $C^1$ boundaries), or cones (for boundary regularity in Lipschitz domains).

By classical results we know that in bounded domains we can achieve exact control,
and therefore any domain satisfying (BU) has to be unbounded. Classical
backward uniqueness results for parabolic equations imply that $\Om=\R^n$
satisfies (BU). It turns out that the half-space $\Om=\R^n_+$ also satisfies
(BU), although this is harder to prove, see \cite{ESS}.  In general, the smaller the domain,
the harder it will be to show that it satisfies (BU). It is immediate that
if $\Om_1\subset\Om_2$ and $\Om_1$ satisfies (BU), then also $\Om_2$
satisfies (BU).

In this paper we consider the question for cones with opening angle $\theta$.
In suitable coordinates
\begin{equation} \label{angledef}
\mathcal O_\theta = \{x=(x_1,x'), x'\in \R^{n-1}, x_1> |x|\cos(\theta/2)\}.
\end{equation}

Luis Escauriaza~\cite{E} recently showed that - surprisingly -  (BU) fails when $\theta<\pi/2$. 
We briefly recall the counterexample. Let us denote by $\Gamma$ the
standard heat kernel, i.\ e.\ $\Gamma(x,t)=(4\pi t)^{-n/2} \exp{(-|x|^2/4t)}$, and recall Appell's transformation
\begin{equation}
\label{appell}
u(x,t)=\Gamma(x,t)v(y,s),\qquad y={x\over t},\quad s={1\over t}\,\,.
\end{equation}

This transformation takes the solutions $u(x,t)$ of the heat equation into the solutions
$v(y,s)$ of the backward heat equation

\begin{equation}
\label{backwardheat}
v_s+\Delta v = 0\,\,.
\end{equation}
By taking $u(x,t)=h(x)$ for a suitable harmonic function $h$ in $\OOa$, we can get a counterexample
to the backward heat equation form of (BU) for $\theta<\pi/2$. In dimension 2 we can use
the real part of the holomorphic function $z\to \exp(-Az^\alpha)$ (for suitable $A>0$ and a parameter $\alpha>2$, $z=x_1+ix_2$) to obtain
an explicit formula:
\begin{equation}
\label{example}
v(y,s)= \Re {1\over s}\, \exp({ -A{{(y_1+iy_2)^\alpha}\over s^\alpha}+{{|y|^2}\over{4s}}})\,\,.
\end{equation}

The function $v(y_1,y_2,s)$ is  bounded in a sector that $|\arctan\frac{y_2}{y_1}|<\pi/(2\alpha)$, away from the origin. We can shift
it to $v(y_1+1,y_2,s)$. The resulted function is bounded in a sector with angle $\pi/\alpha$ satisfying the backward heat equation $v_s+\triangle v=0$, and $v(\cdot, \cdot, 0)=0$. 

We note that it is enough to construct a counterexample in dimension $n=2$. The higher-dimensional
example can then be constructed by simply considering the two-dimensional function
as a function of $n$ variables, independent of $x_3,\dots, x_n$.

Escauriaza's example shows that (BU) fails for $\theta<\pi/2$. Since
we also know that (BU) is true for $\theta=\pi$, it is easy to see that there exists
a borderline angle $\theta_0\in[\pi/2, \pi]$ such that (BU) is true for $\theta>\theta_0$,
and (BU) fails for $\theta<\theta_0$. The borderline case $\theta=\theta_0$ might
perhaps present an extra difficulty.

The main result of this paper is the following:
   \begin{theorem}
   \label{maintheorem}
   The cones $\oth$ satisfy (BU) for 
   $$\theta>2\arccos(1/\sqrt{3})\sim 109.52^{\circ}\,.$$ 
   In other words,
   the critical angle $\theta_0$ introduced above satisfies
    $$\theta_0\le 2\arccos(1/\sqrt{3}).$$
   \end{theorem}

It is tempting to conjecture that $\theta_0=\pi/2$. This is supported by the fact that
$\theta=\pi/2$ is the borderline case for the above construction of Escauriaza,
as can be seen from the Phragm\'en-Lind\"olef principle.

For the classical heat equation, corresponding to the case $b=0$ and $c=0$ in (\ref{heatIntro}),
and $\theta=\pi$ (the half-space), the statement (BU) can be proved by a relatively simple
application of Fourier transformation and some classical complex analysis results,
see \cite{SS1}. We were not able to find such simple proof of the case $b=0,\, c=0$ 
when $\theta<\pi$.

While completing this paper, we learned about reference \cite{GM}.\footnote{We thank Gregory Seregin for pointing out
this article.} Theorem 6 in \cite{GM}  states that for the classical heat equation (corresponding to $b=0,\, c=0$),
 (BU) holds if and only if $\theta\geq \pi/2$. Unfortunately, it seems that the proof is  not available in print.

Our proof of Theorem \ref{maintheorem} relies on two Carleman-type inequalities,
along lines similar to \cite{ESS}. The first inequality, Proposition~\ref{Carleman1},
is taken from \cite{ESS} and is applied in the same way to obtain fast decay rates
for the solutions, see Lemma~\ref{lemmaball}. We note that Carleman inequalities 
of this type can be found already
in~\cite{E2, EV, EF, F}. 

The second inequality, Proposition~\ref{Carleman2}, is the main new tool used in our proof.
The heuristics behind this inequality is somewhat similar to the heuristics behind
the second Carleman-type inequality in \cite{ESS} (Proposition 6.2). However, the proof of
Proposition~\ref{Carleman2} requires a new idea, as for $\theta<\pi$ certain critical terms
appearing in the proofs of the inequalities lose convexity.

In addition to determining the critical angle, another interesting open problem is 
to optimize the assumptions on the coefficients $b(x,t)$ and $c(x,t)$, in the spirit
of~\cite{Koch-Tataru}. For example, it is conceivable that the result remains true
for $b\in L^{n+2}_{x,t}$ and $c\in L^{(n+2)/2}_{x,t}$, but it might be a difficult
problem to decide whether this is the case. 

In what follows we will work with the inequality
\be\label{heatIneqForward}
|u_t-\dd u|\le c_1(|\nn u|+|u|)
\ee
rather than~\rf{heatIntro}. It is not hard to see that when assuming the boundedness of $b(x, t)$ and $c(x, t)$,
the two formulations are equivalent.

\section{Backward uniqueness}
Without loss of generality we assume $T=1$ and work with the backward form of~\rf{heatIneqForward}.
Recall that
\begin{equation}
\mathcal O_\theta = \{x=(x_1,x'), x'\in \R^{n-1}, x_1> |x|\cos(\theta/2)\}.
\end{equation}
Suppose that $u(x, t)$ is a solution to the backward heat equation for some $\theta>2\arccos(1/\sqrt{3})$.
\begin{eqnarray}
|u_t+\Delta u| \leq c_1(|\nabla u|+|u|)\qquad \mbox{in } \OOa\times(0, 1),\label{heat}\\
u(\cdot, 0)=0 \qquad\qquad\quad \mbox{in } \OOa \label{zero}.
\end{eqnarray}
In addition,
\begin{equation}\label{bound}
|u|<M \qquad \mbox{in } \OOa\times(0, 1).
\end{equation}
We will show that  $u\equiv 0$. 

For convenience we will also extend $u(x, t)$ by 0 for $t\leq 0$.

We firstly need the following Carleman inequality from \cite{ESS} (Proposition 6.1), by which we obtain a decay result for solutions of backward heat equation.
\begin{proposition}[\cite{ESS}]
 \label{Carleman1}
  For any function $u\in C^\infty_0(\R^n\times(0, 2); \R^n)$ and any positive number $a$,
	\begin{multline}\label{CarlemanI}
	\int_{\R^n\times(0, 2)}h^{-2a}(t)e^{-\frac{|x|^2}{4t}}\left( \frac{a}{t}|u|^2+|\nabla u|^2\right)dxdt \\ \leq c_0\int_{\R^n\times(0, 2)}h^{-2a}(t)
	e^{-\frac{|x|^2}{4t}}|\partial_t u+\Delta u|^2dxdt,
	\end{multline}
where $c_0$ is a positive absolute constant and $h(t)=te^{\frac{1-t}{3}}$.
\end{proposition}

There are stronger versions of this inequality, see for example \cite{ESS2} inequality (9), and, in particular, \cite{Koch-Tataru}, where inequalities of this type are analyzed in depth.

Lemma \ref{lemmaball} below immediately  implies  exponential decay of the solution $u$ satisfying (\ref{heat})-(\ref{bound}).  The proof is by using the Carleman inequality in Proposition \ref{Carleman1}.  The decay of $u$ enables us to apply the Carleman inequality in Proposition~\ref{Carleman2} and reach the conclusion in Theorem \ref{maintheorem}.

\begin{lemma}\label{lemmaball}
Let  $B_R$ denote the ball with radius $R$ in $\R^n$.  Assume that $R>2$. Consider a function $u$ satisfying the following conditions, with some positive constants $c_1$ and $M$.
\begin{eqnarray}
 |u_t+\triangle u|\leq c_1 (|\nabla u|+|u|) \qquad&&\mbox{ in } B_R\times(0, T),\\
 u(x, 0)=0 \qquad&&\mbox{ in } B_R, \\
 |u|< M \qquad && \mbox { in } B_R\times(0, T).
\end{eqnarray}
Then there exists some constants $\beta$,  $\gamma$, such that for $t\in (0, \gamma)$,
\begin{equation}\label{balldecay}
|u(0, t)|\leq \frac{c_2}{\min\{1, T\}} M e^{-\beta\frac{ R^2}{t}},
\end{equation}
where $\beta$ is a small enough absolute constant,  $c_2$ depends on $c_1$, $\gamma$ depends on $c_1$ and $T$.
 \end{lemma}

Notice that the constants $\beta$, $\gamma$ and $c_2$ do not depend on $R$. It follows immediately the exponential decay of the bounded solution to the equation (\ref{heat}) with respect to the distance to the lateral boundary. We will give the proof of the lemma in the next section.

The second Carleman inequality in Proposition \ref{Carleman2} is a key tool used in our proof of the backward uniqueness in  cones.
We define  the set
$$Q_{\theta}=(\OOa\cap\{x_1>1\})\times(0, 1).$$ The purpose of  ``cutting off the corner'' by requiring $x>1$ is to avoid singularities at the origin.
	\begin{proposition}\label{Carleman2}
	 Let $\phi(x, t)=a\Lambda(t)\varphi(x)+t^2$, where $\Lambda(t)=\dfrac{1-t}{t^{\alpha/2}}$,  and $\varphi(x)=x_1^\alpha-\ep^\alpha r^\alpha$, where $r=|x|$, and $\ep=\cos(\theta/2)$.  For any
	$\ep\in(0, 1/\sqrt{3})$, that is, $\theta\in (2\arccos(1/\sqrt{3}), \pi)$, there exists some $\alpha=\alpha(\ep)\in (1, 2)$ such that the following
	 inequality holds for $u\in C_0^\infty(Q_\theta)$ and $a>a_0$ for some constant $a_0$.
	\begin{eqnarray}\label{Carleman}
	\int_{Q_\theta}e^{2\phi(x, t)}\left[a\left(\Lambda(t)+ \varphi (x) \right) u^2+|\nabla u|^2\right]dxdt\nonumber\\
	\leq 4\int_{Q_{\theta}}e^{2\phi(x, t)}|\partial_t u+\Delta u|^2 dxdt.
	\end{eqnarray}
	\end{proposition}

We apply this Carleman inequality to prove the main result of the paper in the remaining part of this section. The proof of Proposition \ref{Carleman2} is postponed to Section~\ref{CarlemanIIproof}.

For $x\in \OOa$ we denote by $d_\theta(x)$  the distance between $x$ and the boundary of $\OOa$, explicitly given by
	\begin{equation}\label{distance}
	d_\theta(x)=x_1\sin (\theta/ 2) -|x'|\cos(\theta/2).
	\end{equation}
Let $\OOa^{+2}=\{x\in \OOa\, |\,  d_\theta(x)>2\}$. With any other number $c$, the set $\OOa^{+c}$
is defined in the same way.

The next lemma is a consequence of the decay result from Lemma \ref{lemmaball} and Proposition \ref{Carleman2}. It implies Theorem \ref{maintheorem} immediately.

	\begin{lemma}\label{lemmaunique}
	Assume that for some $\theta\in (2\arccos(1/\sqrt{3}), \pi)$ a function $u$ satisfies (\ref{heat}) -- (\ref{bound}), then  there is a number
	$\gamma_1(c_1)$ such that
		\begin{equation}
		u(x, t)\equiv 0
		\end{equation}
	 in $\OOa\times (0,\gamma_1)$.
	\end{lemma}

\begin{proof}
Lemma \ref{lemmaball} implies that
\begin{equation}
	|u(x,t)|\leq c_2M e^{-\beta\frac{d_\theta^2(x)}{t}}
	\end{equation}
for all $(x,t)\in \OOa^{+2} \times(0, \gamma)$.
By local gradient estimates for the heat equation \cite{Lieberman} we can assume that
	\begin{equation}
	|u(x,t)|+|\nabla u(x,t)|\leq c_3M e^{-\frac{\beta d_\theta^2(x)}{2t}}
	\end{equation}
for all $(x, t)\in \OOa^{+3}\times(0, \gamma/2]$.

 By scaling we define a function $v$ by
	 \begin{equation}
	 v(y, s)=u(\lambda y, \lambda^2 s-\gamma_1)
	 \end{equation}
for $(y, s)\in \OOa\times(0, 1)$ with $\lambda=\sqrt{2\gamma_1}$. This function satisfies the relations
	\begin{eqnarray}
	|\partial_s v+\Delta v |\leq c_1\lambda(|\nabla v|+|v|) \qquad \mbox{ in } \OOa\times (0, 1)\label{vlambda},\\
	v(y, s)=0 \qquad \mbox{in } \OOa\times(0, 1/2),
	\end{eqnarray}
and
	\begin{equation}\label{vdecay}
	|v(y, s)|+|\nabla v(y, s)|\leq c_3Me^{-\frac{\beta \lambda^2d_\theta^2(y)}{2(\lambda^2s-\gamma_1)}}\leq c_3Me^{-\beta\frac{d_\theta^2(y)}{2s}}
	\end{equation}
for $1/2<s<1$ and $y\in \OOa^{+3/\lambda}=\{y\in \OOa\, |\,  d_\theta(y)>3/\lambda\}$.

To apply Proposition \ref{Carleman2}, we need certain decay of $|v(y, s)|$ when $|y|$ is large. Notice that the preferred decay can be obtained by considering a cone with slightly smaller opening. Proposition \ref{Carleman2} holds for angles in the interval $(2\arccos(1/\sqrt{3}), \pi)$. We thus consider the median of $\theta$ and $2\arccos(1/\sqrt{3})$,
	$$
	\delta=\frac{\theta+2\arccos(1/\sqrt{3})}{2}.
	$$
In  the smaller cone $\OOd=\{x\in \R^n, x_1>|x|\cos(\delta/2)\}$ we have the estimate
	$
	d_\theta(y)\geq|y|\sin\left(\frac{\theta-\delta}{2}\right)
	$.
It follows (\ref{vdecay}) that
	\begin{equation} \label{vdecay2}
	|v(y, s)|+|\nabla v(y, s)|\leq c_3Me^{-\beta'\frac{|y|^2}{s}}
	\end{equation}
for $1/2<s<1$ and $y\in \OOd\cap\OOa^{+3/\lambda}$,   with the constant $\beta'=\beta\sin^2(\frac{\theta-\delta}{2})/2$.
We can further have
	\begin{equation} \label{vdecay3}
	|v(y, s)|+|\nabla v(y, s)|\leq c'_3Me^{-\beta'\frac{|y|^2}{s}}
	\end{equation}
for $1/2<s<1$ and $y\in \OOd\cap\{y_1> 3/\lambda\}$ for some other constant $c'_3$.

Next we work on the smaller cone $\OOd$ with opening $\delta$, where we have exponential decay (\ref{vdecay3}) and the following properties inherited from $\OOa$.
	\begin{eqnarray}
	|\partial_s v+\Delta v |\leq c_1\lambda(|\nabla v|+|v|) \qquad \mbox{ in } \OOd\times (0, 1),\label{vlambda2}\\
	v(y, s)=0 \qquad \mbox{in } \OOd\times(0, 1/2).\label{vzero}
	\end{eqnarray}

Proposition \ref{Carleman2} requires support condition for the Carleman inequality. For that purpose, let us fix two smooth cut-off functions such that
	$$
	\psi_1(y_1)=\left\{
	\begin{array}{cl}
	0, & y_1<3/\lambda+2,\\
	1, & y_1>3/\lambda+3,
	\end{array}\right.
	$$
		$$
		\psi_2(\tau)=\left\{
		\begin{array}{cl}
		0, &\tau<-3/4,\\
		1, &\tau>-1/2.
		\end{array}\right.
		$$
We set (for the definition of $\phi$, see Proposition~\ref{Carleman2})
	$$
	\phi_B(y, s)=\frac{1}{a}\phi(y, s)-B=(1-s)\frac{y_1^\alpha-\ep^{\alpha}|y|^{\alpha}}{s^{\alpha/2}}+\frac{s^2}{a}-B,
	$$
where $\ep=\cos(\delta/2)$,  $B=\frac{2}{a}\phi(y_\lambda, \frac{1}{2})$, with $y_\lambda=(3/\lambda+3, 0, \cdots,0)$ and
	$$
	\eta(y, s)=\psi_1(y_1)\psi_2(\frac{\phi_B}{B}), \qquad w(y, s)=\eta(y, s)v(y, s).
	$$
The function $w$ is not compactly supported in $Q_{\delta}=(\OOd\cap\{x_1>1\})\times(0, 1)$. However,  it follows from (\ref{vdecay}) and the special structure of the weight function $\phi$ in Proposition~\ref{Carleman2} that, with $w$ replacing $u$ in (\ref{Carleman}),  integrals on both sides  converge. If we multiply $w$ by an additional cut-off function  $\xi$ such that
	          $$
		\xi(x)=\left\{
		\begin{array}{cl}
		1, &|x|< R\\
		0, &|x|> 2R
		\end{array}\right.
		$$
and $|\nabla\xi|< c/R$, $|\nabla^2 \xi|<c/R^2$, then apply Proposition~\ref{Carleman2} to the compact supported function $w\xi$,  and let $R\to\infty$, we finally obtain
	\begin{multline}
	\int_{Q_{\delta}}e^{2a\phi_B}\left[a\left(\Lambda(s)+ \varphi \right) w^2+|\nabla w|^2\right] dyds\\
	\leq 4\int_{Q_{\delta}}e^{2a\phi_B}|\partial_s w+\Delta w|^2 dyds.
	\end{multline}
From	 (\ref{vlambda2}) we have
	\begin{eqnarray}
	|\partial_s w+\triangle w|&\leq& c_1\lambda(|\nabla w|+|w|)\nonumber\\
	&&+c_{4}(|\nabla v|+|v|) (|\partial_s \eta+|\nabla\eta|+|\Delta\eta|).
	\end{eqnarray}
For large enough parameter $a$ we have that $a(\Lambda(s)+\varphi(y))>1$ in the support of $w$, where $\phi_B/B\geq -3/4$ because of the definition of $\psi_2$. In addition, we take $\gamma_1(c_1)$ small enough such that $16c_1^2\lambda^2<1/2$. We then have
	\begin{eqnarray}
	I&\equiv&\int_{Q_\delta}e^{2a\phi_B}\left(w^2+|\nabla w|^2\right)dyds\\
	&\leq&32c^2_{4} \int_{Q_\delta}e^{2a\phi_B}(|v|^2+|\nabla v|^2)(|\partial_s \eta|+|\nabla\eta|+|\Delta\eta|)^2dyds.\label{I_unique}
	\end{eqnarray}

To estimate the right hand side,  we need to look into details of derivatives of the function $\eta(y, s)$.
In view of the definitions of $\psi_1$ and $\psi_2$,  the support of derivatives of $\eta(y, s)$ is the closure of the set
	\begin{multline*}
	\{ y_1> 3/\lambda+2, \, -3B/4<\phi_B<-B/2 \}\\
	\cup\{3/\lambda+2<y_1<3/\lambda+3, \, \phi_B>-B/2  \}.
	\end{multline*}
However, the second set has empty intersection with $\OOd\times(1/2, 1)$, where the function $v$ is nonzero. Hence the support of the term
$(|\nabla v|+|v|) (|\partial_s \eta+|\nabla\eta|+|\Delta\eta|)$ is closure of the set
	$$
	\{  y_1>3/\lambda+2, 1/2<s<1, -3B/4 <\phi_B(y, s)<-B/2  \},
	$$
of which we denote by $\chi(y, s)$ the characteristic function. 

Next we estimate the derivatives of $\eta(y, s)$ in the support of the term $(|\nabla v|+|v|) (|\partial_s \eta+|\nabla\eta|+|\Delta\eta|)$. Recall that $\phi(y, s)=a\Lambda(s)\varphi(y)+s^2$, where  $\Lambda(s)=\dfrac{1-s}{s^{\alpha/2}}$ and $\varphi(y)=y_1^\alpha-\ep^\alpha |y|^\alpha$. The function $\Lambda(s)$  and the derivative $\Lambda'(s)$ are bounded for  $s\in(1/2, 1)$. The function  $\varphi$ and its derivatives up to the second order are bounded by a constant multiple of $|y|^\alpha$.
The cut-off functions $\psi_1$ and $\psi_2$ and derivatives up to the second order  are  bounded by some absolute constant. The value of $B$ is bounded from below regardless of the value of the parameter $a$. Thus
	\begin{equation}\label{eta}
	(|\partial_s \eta|+|\nabla\eta|+|\Delta\eta|)^2<c_{5}|y|^{2\alpha}.
	\end{equation}

We now estimate (\ref{I_unique}) by using (\ref{vdecay2}) and (\ref{eta}). We see that
	$$
	I\leq c_{6} Me^{-Ba} \int_{Q_{\delta}}|y|^{2\alpha} e^{-2\beta' \frac{|y|^2}{s}}\chi(y, s)dyds
	$$
for some constant $c_{6}$.
The last integral is bounded. Passing to the limit as $a\to\infty$ we see that $v(y, s)=0$ where $\phi_B(y, s)>0$ and $1/2<s<1$. Using the property of unique continuation across the spatial boundaries (see Theorem 4.1 in \cite{ESS}), we show that $v(y, s)=0$ if $y\in \OOa$ and $0<s<1$. This proves the lemma.
\end{proof}

\section{Proof of Lemma 2.2}

The proof of Lemma~\ref{lemmaball} is based on the Carleman inequality in \cite{ESS}, which we quoted in Proposition~\ref{Carleman1}.

\begin{proof} [Proof of Lemma \ref{lemmaball}]The proof is similar to the one in \cite{ESS}. We still include the proof here for the convenience of the reader.

In what follows, we always assume that the function $u$ is extended by zero to negative values of $t$.

The assumption that $R>2$ results in no loss, since the conclusion is only useful when $R$ is large.  According to the local gradient estimates of the heat equation \cite{Lieberman},  in the smaller cylinder $(x, t) \in B_{R-1}\times(0, T/2)$, we can assume that
	\begin{equation}
	|u|+|\nabla u(x, t)|\leq \frac{c_7}{\min\{1, T\}} M
	\end{equation}
with some absolute constant $c_7$.

We fix $t\in (0,\, \min\{1, T\}/12)$ and introduce a new function $v$ by the usual parabolic scaling:
	$$
	v(y, s)=u(\lambda y, \lambda^2s-t/2).
	$$
The function $v$ is well defined on the set $Q_\rho=B_\rho\times(0, 2)$, where $\rho=(R-1)/\lambda$ and $\lambda=\sqrt{3t} \in (0,\,  \min\{1,\sqrt{T}\}/2)$.
We have the following relations for $v$.
	\begin{equation}\label{vheat}
	|\partial_s v+\Delta v|\leq c_1\lambda(|\nabla v|+ |v|),
	\end{equation}
	\begin{equation}\label{vgradv}
	|v(y, s)|+|\nabla v(y, s)| < \frac{c_7}{\min\{1, T\}}M
	\end{equation}
for all $(y, s)\in Q_\rho$,
	\begin{equation}
	v(y, s)=0
	\end{equation}
for $(y, s)\in B_\rho\times(0, 1/6]$.

By the assumption that $R>2$ and $\lambda <1/2$, we have $\rho>2$. In order to apply Proposition \ref{Carleman1}, we take two smooth cut-off functions in $Q_\rho$:
	$$
	\psi_\rho(y)=\left\{
	\begin{array}{cl}
	0, &|y|>\rho-1/2,\\
	1, & |y|<\rho-1,
	\end{array}\right.
	$$
		$$
		\psi_t(s)=\left\{
		\begin{array}{cl}
		0, &7/4<s<2,\\
		1, & 0<s<3/2.
		\end{array}\right.
		$$
By assumption, these functions take values in $[0, 1]$ and are such that $|\nabla^k \psi_\rho|<C_k$, $k=1, 2$,  and $|\partial_s \psi_t| < C_0$. We set $\eta(y, s)=\psi_\rho(y)\psi_t(s)$ and
	\begin{equation}
	w(y, s)=\eta(y, s) v(y, s).
	\end{equation}
It follows from (\ref{vheat}) that
	\begin{equation}\label{eqw}
	|\partial_s w+\Delta w|\leq c_1\lambda(|\nabla w|+|w|)+ c_8 \chi (|\nabla v|+|v|),
	\end{equation}
where $c_8$ is a positive constant depending only on $c_1$ and $C_k$, $k=0, 1, 2$, and $\chi(y, s)=1$ for $(y, s)\in \omega=\{\rho-1<|y|<\rho, 0<s<2\}\cup\{|y|<\rho-1, 3/2<s<2\}$, and $\chi(y, s)=0$ for $(y, s)\notin \omega$. The set $\omega$ is where the cut-off function $\eta$ is not constantly $1$ in $Q_\rho$. Obviously, the function $w$ is compactly supported on $\R^2\times(0, 2)$. The inequality (\ref{CarlemanI}) also holds for scale valued functions.  Therefore we may apply Proposition \ref{Carleman1} and obtain
	\begin{multline}
	\int_{Q_\rho}h^{-2a}(s)e^{-\frac{|y|^2}{4s}}\left( \frac{a}{s}|w|^2+|\nabla w|^2\right)dyds \\ \leq c_0\int_{Q_\rho}h^{-2a}(s)
	e^{-\frac{|y|^2}{4s}}|\partial_s w+\Delta w|^2dyds.
	\end{multline}
Taking $a>2$, and applying (\ref{eqw})
we finally obtain that
	\begin{equation}\label{smallambda}
	I\equiv\int_{Q_\rho} h^{-2a}(s) e^{-\frac{|y|^2}{4s}}(|w|^2+|\nabla w|^2) dyds\leq 4c_0(c_1^2\lambda^2I+c_8^2I_1),
	\end{equation}
where
	$$
	I_1=\int_{Q_\rho}\chi(y, s)h^{-2a}(s)e^{-\frac{|y|^2}{4s}}(|\nabla v|^2+|v|^2)dyds.
	$$
Taking a sufficiently small value for $\gamma=\gamma(c_1)$ such that in the range $\lambda\in (0, \gamma)$, we can assume that the inequality $4c_0c_1^2\lambda^2\leq 1/2$ holds, and therefore (\ref{smallambda}) implies that
	\begin{equation}
	I\leq 8c_0c_8^2 I_1.
	\end{equation}
Notice that near the origin $\{y=0, s=0\}$, where the parametric function $h^{-2a}(s)e^{-\frac{|y|^2}{4s}}$ is not integrable,  the characteristic function $\chi$ is $0$.
By (\ref{vgradv}) we have
	\begin{eqnarray}
	I_1&\leq&\frac{c_7^2M^2}{\min\{1, T^2\}}\left\{\int_{3/2}^{2}\int_{|y|<\rho-1}h^{-2a}(s)e^{-\frac{|y|^2}{4s}}dyds\right. \nonumber\\
	      &  &\left.+\int_0^{2}\int_{\rho-1<|y|<\rho}h^{-2a}(s)e^{-\frac{|y|^2}{4s}}dyds\right\}\\
	      &\leq&\frac{c_9M^2}{\min\{1, T^2\}}\left[h^{-2a}(3/2)+\int_{0}^{2}h^{-2a}(s)e^{-\frac{(\rho-1)^2}{4s}}ds\right].\label{I_1}
	\end{eqnarray}

Using (\ref{I_1}) we obtain the estimate
	\begin{eqnarray*}
	D&\equiv&\int_{B_1}\int_{1/2}^1|w|^2 dyds=\int_{B_1}\int_{1/2}^1|v|^2dyds \\
	&\leq& c_{10}\int_{Q_\rho} h^{-2a}(s) e^{-\frac{|y|^2}{4s}}(|w|^2+|\nabla w|^2) dyds\\
	&\leq& \frac{c_{11}M^2}{\min\{1, T^2\}}\left[h^{-2a}(3/2)+\int_{0}^{2}h^{-2a}(s)e^{-\frac{\rho^2}{16s}}\, ds\right]\\
	&=& \frac{c_{11}M^2}{\min\{1, T^2\}} e^{-\beta\rho^2}\left[h^{-2a}(3/2)e^{\beta\rho^2}+\int_{0}^{2}h^{-2a}(s)e^{\beta\rho^2-\frac{\rho^2}{16s}}\,ds\right].
	\end{eqnarray*}
We take $\beta<1/64$ and then let
	$$
	a=\beta\rho^2/(2\log h(3/2)).
	$$
This choice of $a$ leads to the estimate
	$$
	D\leq \frac{c_{11}M^2}{\min\{1, T^2\}}e^{-\beta\rho^2}\left[1+\int_{0}^{2}g(s)ds\right],
	$$
where $g(s)=h^{-2a}(s)e^{-\frac{\rho^2}{32s}}$. By simple calculation we have that
	$$
	g'(s)=h^{-2a}(s)e^{-\frac{\rho^2}{32s}}\left[-\frac{\beta\rho^2}{\log h(3/2)}\left(\frac{1}{s}-\frac{1}{3}\right)+\frac{\rho^2}{32s^2}\right].
	$$
One can readily verify that $g(2)<1$ and $g'(s)\geq 0$ for any $s\in (0, 2)$ if $\beta<\frac{1}{64}\log h(3/2)$. Therefore,
	$$
	D\leq 3 \frac{c_{11}M^2}{\min\{1, T^2\}}e^{-\beta \rho^2}=3 \frac{c_{11}M^2}{\min\{1, T^2\}}e^{-\beta\frac{ R^2}{12t}}.
	$$
On the other hand, the regularity theory implies that
	$$
	|u(0, t)|^2=|v(0, 1/2)|^2\leq c_{12} D.
	$$
Finally we obtain that $$|u(0, t)|\leq\frac{c_2}{\min\{1, T\}}M e^{-\beta\frac{ R^2}{24t}}.$$
\end{proof}


\section{Proof of Proposition 2.3}\label{CarlemanIIproof}

One of the difficulties in the proof of the Carleman inequality in Proposition~\ref{Carleman2} is that 
-- by comparison with the case $\theta=\pi$ -- 
some loss of the convexity of the weight $\vf$  cannot be avoided. Therefore we have to investigate 
in more detail some of the terms which can be neglected when $\theta\ge\pi$.  

\begin{proof}[Proof of Proposition~\ref{Carleman2}]
 Let $u$ be an arbitrary function in $C_0^\infty(Q_{\theta})$ and $v=e^{\phi}u$. Then
\begin{equation}Lv\equiv e^\phi(\partial_t u+\Delta u)=\Delta v+|\nabla\phi|^2v-\partial_t\phi v + \partial_t v -2\nabla\phi\nabla v-\Delta\phi v.\end{equation}

We decompose $L$ into symmetric and skew-symmetric parts $$L=S+A,$$ where
	\begin{equation}
	Sv=\Delta v+|\nabla\phi|^2v-\partial_t\phi v
	\end{equation}
and
	\begin{equation}
	Av=\partial_t v -2\nabla\phi\nabla v-\Delta\phi v.
	\end{equation}

The right hand side of the inequality (\ref{Carleman}) is
	\begin{equation}
	 \int |Lv|^2 dxdt=\int |Sv|^2 dxdt+\int |Av|^2 dxdt+\int ([S, A]v) v dxdt,
	\end{equation}
where $[S, A]=SA-AS$ is the commutator of $S$ and $A$. By simple calculations we have that
	\begin{eqnarray}\label{commutator}
	([S, A]v, v)&=&\int 4\phi_{,kl}v_{,k}v_{,l}dxdt\label{commutatorDerivative}\\
	&+&\int\left( 2\nabla\phi \nabla|\nabla\phi|^2-\Delta^2\phi+\partial^2_t\phi-2\partial_t|\nabla\phi|^2\right)|v|^2 dxdt.
	\end{eqnarray}

The Hessian of the function $\phi=a\Lambda(t)(x_1^\alpha-\ep^\alpha r^\alpha)$ is not positive-definite. The loss of convexity is not avoidable with an angle smaller than $\pi$. To compensate the term $\int 4\phi_{,kl}v_{,k}v_{,l}dxdt$ in (\ref{commutator}) we introduce a function $F(x, t)$ to be  determined.
	\begin{eqnarray*}
	(Sv, Fv)&=&\int \Delta v F v+(|\nabla \phi|^2 -\partial_t \phi)Fv^2 dx\\
	&=&\int -F|\nabla v|^2+(\frac{1}{2}\triangle F+|\nabla \phi|^2 F-\partial_t \phi F)v^2 dxdt.
	\end{eqnarray*}
Cauchy-Schwartz inequality implies that
	\begin{eqnarray}\label{symm}
	(Sv, Sv)&\geq&-(Sv, Fv)-\frac{1}{4}\int F^2v^2 dxdt\nonumber\\
	&\geq& \int F|\nabla v|^2-(\frac{1}{2}\triangle F+|\nabla \phi|^2 F-\partial_t \phi F+\frac{1}{4}F^2)v^2 dxdt.
	\end{eqnarray}
Combining (\ref{commutator}) and (\ref{symm}) we have
	\begin{eqnarray}
	\lefteqn{([S, A]v, v)+(Sv, Sv)\geq
	 \int 4\phi_{,kl}v_{,k}v_{,l}+ F |\nabla v|^2 dxdt} \label{gradv}\\
	&&+\int\left( 2\nabla\phi \nabla|\nabla\phi|^2-\Delta^2\phi+\partial^2_t\phi-2\partial_t|\nabla\phi|^2\right)v^2 dxdt\\
	&& +\int-(\frac{1}{2}\triangle F+|\nabla \phi|^2 F-\partial_t \phi F+\frac{1}{4}F^2)v^2 dxdt.
	\end{eqnarray}
By calculation the Hessian of $\varphi$ is
	\begin{eqnarray}
	D^2\varphi(x)
	&=& \alpha(\alpha-1) \left(\begin{array}{cccc}
	x_1^{\alpha-2}
	&0&\cdots&0\\
	0&0&\cdots&0\\
	\vdots& \vdots& \ddots&\vdots\\
	0&0&\cdots&0
	\end{array}
	\right)
	-\alpha\ep^\alpha r^{\alpha-2}E_n\\
	&&+\alpha(2-\alpha)\ep^\alpha r^{\alpha-4}x^{T}x,
	\end{eqnarray}
where $E_n$ reprensents the $n$ dimensional identity matrix, $x=(x_1, \dots, x_n)$ is the row vector and $x^{T}$ denotes the transpose of $x$.
It is easy to see that
          $$
	D^2\varphi(x)+\alpha \ep^{\alpha}r^{\alpha-2}E_n\geq 0.
	$$
We thus let  $f(x)=\alpha\ep^{\alpha} r^{\alpha-2}$ and
	\begin{equation}
	F(x, t)=4a\Lambda(t)f(x)+1.
	\end{equation}
With this choice of $F(x, t)$, the right hand side of line (\ref{gradv}) is positive and
	$$
	 \int 4\phi_{,kl}v_{,k}v_{,l}+ F |\nabla v|^2 dxdt \geq \int |\nabla v|^2dxdt.
	$$
Grouping the remaining terms according to the orders of the parameter $a$, with $A_3$ denoting the terms with $a^3$ and etc.,  we have that
	\begin{eqnarray}\label{commsymm}
	([S, A]v, v)+(Sv, Sv) &\geq& \int |\nabla v|^2 dxdt\nonumber\\
	&&+\int (A_3+A_2+A_1+A_0)v^2 dxdt,
	\end{eqnarray}
where
	\begin{eqnarray}
	A_3&=&2\nabla\phi \nabla|\nabla\phi|^2-4a^3\Lambda^3 f|\nabla\varphi|^2,\\
	A_2&=&- 4a^2\Lambda\Lambda' |\nabla\varphi|^2+4a^2\Lambda\Lambda'\varphi f-4a^2\Lambda^2(t) f^2-a^2 \Lambda^2|\nabla\varphi|^2,\\
	A_1&=&-a\Lambda\Delta^2\varphi+a\Lambda''\varphi-2a\Lambda(t) \Delta f+a \Lambda' \varphi-2a\Lambda(t)f+8a t\Lambda f,\\
	A_0&=&7/4+2t.
	\end{eqnarray}

We analyze $A_3$ first. With $a^3$ as a coefficient $A_3$ must be non-negative in the set $Q_\theta$.  By simple calculation we can see that $2\nabla\phi \nabla|\nabla\phi|^2> 0$ is equivalent to $\ep<1/\sqrt{2}$, or $\theta>\pi/2$. However, to compensate the loss of convexity in (\ref{commutatorDerivative}) we introduced a negative term $-4a^3\Lambda^3 f|\nabla\varphi|^2$.  By letting $x_1/r\to\ep$, we see that  $\ep< 1/\sqrt{3}$ is a necessary condition for $A_3\geq 0$. Next we show that under the condition $\ep<1/\sqrt{3}$, we indeed have that $A_3\geq0$. Denoting by $\nabla\varphi^{T}$ the transpose of the row vector $\nabla\varphi$,  we notice that
	\begin{equation}
	A_3=
	4a^3\Lambda^3(t)
	\nabla\varphi
	\left(D^2\varphi(x)-\alpha\ep^{\alpha}r^{\alpha-2}\right)
	\nabla\varphi^{T}.
	\end{equation}
It is easy to see that
	\begin{eqnarray}
	\lefteqn{D^2\varphi(x)-\alpha\ep^{\alpha}r^{\alpha-2}
	}&&\nonumber\\
	 &\geq&\alpha(\alpha-1) \left(\begin{array}{cccc}
	x_1^{\alpha-2}
	&0&\cdots&0\\
	0&0&\cdots&0\\
	\vdots& \vdots& \ddots&\vdots\\
	0&0&\cdots&0
	\end{array}
	\right)
	-2\alpha\ep^\alpha r^{\alpha-2}E_n.
	\end{eqnarray}	
By the fact that  $x_1^{\alpha-2}>r^{\alpha-2} $ (since $\alpha<2$), we have
\begin{equation}
	D^2\varphi(x)-\alpha\ep^{\alpha}r^{\alpha-2}
	\geq
	\alpha r^{\alpha-2}
	\left(\begin{array}{cc}
	\alpha-1-2\ep^{\alpha} &0\\
	0 & -2\ep^{\alpha}E_{n-1}
	\end{array}\right),
\end{equation}
where $E_{n-1}$ is the $n-1$ dimensional identity matrix.

The first derivatives of $\varphi$ are as follows.
	$$\varphi_{,1}(x,t)=\alpha x_1^{\alpha-1}-\alpha\ep^\alpha r^{\alpha-2}x_1, \quad \varphi_{,k}(x,t)=-\alpha\ep^\alpha r^{\alpha-2}x_k, \quad k=2, \dots, n.$$
We notice that $\varphi_{,1}\geq\alpha (1-\ep^{\alpha}) x_1^{\alpha-1}$. Thus
	\begin{equation*}A_3\geq4a^3\Lambda^3(t)\alpha^3 r^{\alpha-2}\left[(\alpha-1- 2\ep^{\alpha})(1-\ep^{\alpha})^{2}x_1^{2\alpha-2}-
	2\ep^{3\alpha}r^{2\alpha-4}|x'|^2
	\right].
	\end{equation*}
Taking into account that $x_1/r>\ep$ and $|x'|^2/r^2<1-\ep^2$,
	\begin{equation*}
	A_3\geq4a^3\Lambda^3(t)\alpha^3 r^{3\alpha-4} \ep^{2\alpha-2}\left[
	(\alpha-1- 2\ep^{\alpha})(1-\ep^{\alpha})^{2}
	- 2\ep^{\alpha+2}(1-\ep^2)\right].
	\end{equation*}
Let us denote the quantity in the bracket above by $m(\alpha, \ep)$.
	$$
	m(\alpha, \ep)=(\alpha-1- 2\ep^{\alpha})(1-\ep^{\alpha})^{2}
		- 2\ep^{\alpha+2}(1-\ep^2).
	$$
$A_3$ is non-negative if $m(\alpha, \ep)$ is. In the set $(\alpha, \ep)\in(1, 2)\times(0, 1/\sqrt{3})$, the function $m(\alpha, \ep)$ is monotone increasing with respect to $\alpha$ and is monotone decreasing with respect to $\ep$.  Notice that
	$$
	m(2, 1/\sqrt{3})=0.
	$$
Hence for any $\ep$ smaller than and close to $1/\sqrt{3}$,  we can choose a corresponding $\alpha<2$ such that  $m(\alpha, \ep)\geq 0$, and in turn
$A_3\geq 0$. We fix this $\alpha$ in the rest of the proof.

 For the estimate of $A_2$ we notice that $\nabla v=(\nabla\phi )v+e^{\phi}\nabla u$, as $v=e^\phi u$. To bound $|e^{\phi}\nabla u|$, we want to apply the inequality  $|e^{\phi}\nabla u|^2/2\leq|\nabla v|^2+|\nabla\phi|^2 v^2$. We thus look at the inequality (\ref{commsymm}) in the following way.
          \begin{eqnarray*}
	([S, A]v, v)+(Sv, Sv) &\geq& \int (|\nabla v|^2+|\nabla\phi|^2 v^2) dxdt\nonumber\\
	&&+\int [A_3+(A_2-|\nabla\phi|^2)+A_1+A_0]v^2 dxdt.
	\end{eqnarray*}
Next we estimate $A_2-|\nabla\phi|^2$.
	\begin{equation}
	A_2-|\nabla\phi|^2 =-4a^2\Lambda(t)\Lambda'(t)\left( \left( 1+\frac{\Lambda(t)}{2\Lambda'(t)}\right)|\nabla\varphi|^2-\varphi f+
	\frac{\Lambda(t)}{\Lambda'(t)} f^2 \right).
	\end{equation}
$\Lambda(t)=\frac{1-t}{t^{\alpha/2}}$ and $\Lambda'(t)=-\frac{\alpha/2+(1-\alpha/2)t}{t^{\alpha/2+1}}$. $|\Lambda(t)/\Lambda'(t)|<\frac{1}{2\alpha}$.
	\begin{equation}
	A_2-|\nabla\phi|^2 \geq-4a^2\Lambda(t)\Lambda'(t)\left( \left( 1-\frac{1}{4\alpha}\right)|\nabla\varphi|^2-\varphi f-
	\frac{1}{2\alpha} f^2 \right).
	\end{equation}
Notice that $|\nabla\varphi|^2\geq |\varphi_{,1}|^2\geq\alpha^2(1-\ep^{\alpha})^2x_1^{2\alpha-2}$ and $\varphi f\leq \alpha\ep^\alpha(1-\ep^{\alpha})x_1^{2\alpha-2}$. For the choice of $\alpha$ we made above (see the expression of $m(\alpha, \ep)$), $\ep^\alpha\leq(\alpha-1)/2$.  Taking into account that $ r\geq x_1>1$,
	\begin{equation}
	A_2-|\nabla\phi|^2 \geq-4a^2\Lambda(t)\Lambda'(t)\left( \frac{3}{8}x_1^{2\alpha-2}-\frac{1}{4}r^{2\alpha-4}\right)\geq-\frac{a^2}{2}\Lambda(t)\Lambda'(t)x_1^{2\alpha-2}.
	\end{equation}

Finally, we estimate $A_1$. Recall that

        $$
        A_1 =-a\Lambda\Delta^2\varphi+a\Lambda''\varphi-2a\Lambda(t) \Delta f+a \Lambda' \varphi-2a\Lambda(t)f+8a t\Lambda f.
        $$
A simple observation  is that
\begin{equation}
A_1 \geq a(\Lambda''+\Lambda') \varphi -a\Lambda(\Delta^2\varphi+2\Delta f+2f),
\end{equation}
and  $\Lambda''(t)+\Lambda'(t)>\alpha-1>0$. $|-a\Lambda(\Delta^2\varphi+2\Delta f+2f)|< C a\Lambda$,  thus under the control of $A_2$. Consequently,  there exists some constant $a_0$ depending on $\varphi$ such that for $a>a_0$,
	\begin{equation}
	A_2+A_1\geq -\frac{a^2}{4}\Lambda(t)\Lambda'(t)x_1^{2\alpha-2}+a(\alpha-1)\varphi.
	\end{equation}
In addition $|\Lambda'(t)|\geq1$ and $x_1>1$.
We thus have
          \begin{eqnarray*}
	([S, A]v, v)+(Sv, Sv) &\geq& \int (|\nabla v|^2+|\nabla\phi|^2 v^2) dxdt\\
	&&+\int\frac{a^2}{4}\Lambda(t)v^2 dxdt+\int a(\alpha-1)\varphi v^2 dxdt
	\end{eqnarray*}
In turn
	\begin{eqnarray}
	\int_{Q_\theta}e^{2a\Lambda(t)\varphi(x)+2t^2}\left[\left(\frac{a^2}{4}\Lambda(t)+ a(\alpha-1)\varphi \right) u^2+\frac{1}{2}|\nabla u|^2
	\right]dxdt\nonumber\\
	\leq \int_{Q_{\theta}}e^{2a\Lambda(t)\varphi(x)+2t^2}|\partial_t u+\Delta u|^2 dxdt.
	\end{eqnarray}
To simplify, we can assume $\alpha>3/2$, and $a>2$, it follows that
	\begin{eqnarray}
	\int_{Q_\theta}e^{2a\Lambda(t)\varphi(x)+2t^2}\left[a\left(\Lambda+\varphi \right) u^2+|\nabla u|^2\right]dxdt\nonumber\\
	\leq 4\int_{Q_{\theta}}e^{2a\Lambda(t)\varphi(x)+2t^2}|\partial_t u+\Delta u|^2dxdt.
	\end{eqnarray}
\end{proof}

\begin{remark}
Our proof here is quite elementary, in that is does not use micro-local analysis. It is likely that more sophisticated methods, such as those in \cite{Tataru} can be used to obtain additional insights into the inequality and its possible sharper versions.
\end{remark}

{\noindent\bf Acknowledgement}: The authors would like to thank Luis Escauriaza for valuable suggestions and comments .


\end{document}